# An improved car-following model considering variable safety headway distance

Yu-han Jia[a,*], Jian-ping Wu[a]

([a] Department of Civil Engineering, Tsinghua University, Beijing 100084, China)

**Abstract:** Considering high speed following on expressway or highway, an improved car-following model is developed in this paper by introducing variable safety headway distance. Stability analysis of the new model is carried out using the control theory method. Finally, numerical simulations are implemented and the results show good consistency with theoretical study.

**Re´sume´:** A l'égard du phénomène de car-following, c'est-à-dire la suite de voiture, à grande vitesse sur l'autoroute, cet article présente un modèle car-following basé sur la distance de sécurité minimum de variable. Ensuite, l'analyse de la stabilité du nouveau modèle en utilisant la théorie de la commande est effectué, suivie par un grand nombre d'expérimentations simulées mises en oeuvre. Enfin, les résultats montrent que l'étude théorique du nouveau modèle est correcte.

---
[*] Corresponding author: Yu-han Jia
  E-mail address: jiajiayuyuhanhan@163.com

# 1. Introduction

Car-Following theory is a method to describe how vehicles follow one another on a roadway. In the model cars are regarded as discrete and interacting particles without overtaking to analyze traffic flow characteristics. So far, a variety of models have been developed, including safety distance model, stimulus-response model and fuzzy logic based model [1-8]. Among those models, the optimal velocity model (OVM) is well known for accuracy and rationality [3,4,9]. Afterwards the OVM has been extended by introducing generalized force, full velocity difference and lateral effects [10-14]. However, the models mentioned cannot be used to precisely reflect the car-following phenomenon in expressway or highway with higher speed and larger headway (mostly longer than 10m), in which case the hyperbolic tangent part ($\tanh(\cdot)$) in OVM remains as a constant. Therefore the safety headway distance parameter in OVM needs to be modified.

Based on previous work, this paper investigates a new car-following model considering variable safety headway distance (VSHD). In section 2 the new model is developed and stability analysis is carried out in section 3. Then numerical simulation experiments are performed to verify the theoretical study in section 4. The summary is given in section 5.

# 2. Improved model

The typical OVM is presented as

$$\frac{d^2 x_n(t)}{dt^2} = \alpha[V^{op}(\Delta x_n(t)) - v_n(t)] \quad (1)$$

where $x_n(t)$ and $v_n(t)$ are the position and velocity of the nth vehicle; $\Delta x_n(t)$ is the headway distance between the nth and its leading vehicle; $\alpha$ is the sensitivity parameter of the driver; $V^{op}(\cdot)$ is the OV function described as

$$V^{op}(\Delta x_n(t)) = \frac{v_{\max}}{2}[\tanh(\Delta x_n(t) - h_c) + \tanh(h_c)] \quad (2)$$

where $v_{\max}$ is the maximum velocity on a particular roadway; $h_c$ means the safety headway distance.

However, as noticed in the study on expressway where vehicles have higher speed and larger headway, the $\Delta x_n(t) - h_c$ part is larger than on regular road which makes the $\tanh(\cdot)$ remain as 1. Hence the fol-lowing vehicle gets little or no influence despite $\Delta x_n(t)$ is changing, which contradicts the real traffic situation.

To avoid the mentioned problem, the fixed parameter $h_c$ in OV function is replaced by VSHD as $h_f = bv_n(t)t_s + h_c$, where $t_s$ is the time step unit. The physical meaning of $h_f$ is that the acceptable safety headway distance of a driver is dynamic changing based on $v_n(t)$ level.

Referring to previous study [14,15], the velocity difference between the nth and its leading vehicle $\Delta v_n(t)$ is introduced, making the new function as

$$\begin{cases} \dfrac{d^2 x_n(t)}{dt^2} = \alpha[V_{new}^{op}(\Delta x_n(t), v_n(t)) - v_n(t)] \\ \qquad\qquad + \lambda \Delta v_n(t) \\ V_{new}^{op}(\Delta x_n(t), v_n(t)) = \dfrac{v_{\max}}{2}[\tanh(\Delta x_n(t) - h_f) \\ \qquad\qquad + \tanh(h_f)] \\ h_f = bv_n(t)t_s + h_c \end{cases} \quad (3)$$

## 3. Stability analysis

According to stability analysis method [15-17], the stable condition of the modified car-following model can be shown as

$$\begin{cases} \dfrac{dv_n(t)}{dt} = \alpha[V_{new}^{op}(\Delta x_n(t), v_n(t)) - v_n(t)] \\ \qquad\qquad + \lambda \Delta v_n(t) \\ \dfrac{dy_n(t)}{dt} = v_{n+1}(t) - v_n(t) \end{cases} \quad (4)$$

where $y_n(t) = x_{n+1}(t) - x_n(t)$. We assume the first vehicle is not influenced by others and runs constantly at speed $v_0$, then the steady state is given by

$$\left[v_n^*(t), y_n^*(t)\right]^T = \left[v, V_{new}^{op^{-1}}(v_0)\right]^T . \quad (5)$$

The linearized system of (9) can be calculated around steady state (10) as

$$\begin{cases} \dfrac{d\delta v_n(t)}{dt} = \alpha\left[\delta y_n(t)\Lambda_1 + \delta v_n(t)\Lambda_2 - \delta v_n(t)\right] \\ \qquad\qquad + \lambda(\delta v_{n+1}(t) - \delta v_n(t)) \\ \dfrac{d\delta y_n(t)}{dt} = \delta v_{n+1}(t) - \delta v_n(t) \end{cases} \quad (6)$$

where $\delta v_n(t) = v_n(t) - v$, $\delta y_n(t) = y_n(t) - V^{op-1}(v)$, partial derivatives $\Lambda_1 = \dfrac{\partial V_{new}^{op}(\Delta x_n(t), v_n(t))}{\partial \Delta x_n(t)}\Big|_{y_n(t)=V_{new}^{op^{-1}}(v_0)}$, and $\Lambda_2 = \dfrac{\partial V_{new}^{op}(\Delta x_n(t), v_n(t))}{\partial v_n(t)}\Big|_{y_n(t)=V_{new}^{op^{-1}}(v_0)}$.

After Laplace transformation, we have

$$\begin{bmatrix} V_n(s) \\ Y_n(s) \end{bmatrix} = \dfrac{V_{n+1}(s)}{p(s)}\begin{bmatrix} s & \alpha\Lambda_1 \\ -1 & s+\alpha+\lambda-\alpha\Lambda_2 \end{bmatrix}\begin{bmatrix} \lambda \\ 1 \end{bmatrix}, \quad (7)$$

where $V_n(s) = L(\delta v_n(t))$, $Y_n(s) = L(\delta y_n(t))$, $L(\cdot)$ is the Laplace transformation and the characteristic polynomial is $p(s) = s^2 + s\lambda + s\alpha + \alpha\Lambda_1 - s\alpha\Lambda_2$. Then the transfer function can be obtained as

$$G(s) = \begin{bmatrix} 1 & 0 \end{bmatrix} \frac{1}{p(s)} \begin{bmatrix} s & \alpha\Lambda_1 \\ -1 & s+\alpha+\lambda-\alpha\Lambda_2 \end{bmatrix} \begin{bmatrix} \lambda \\ 1 \end{bmatrix}$$
$$= \frac{s\lambda + \alpha\Lambda_1}{p(s)} = \frac{s\lambda + \alpha\Lambda_1}{s^2 + s\lambda + s\alpha - s\alpha\Lambda_2 + \alpha\Lambda_1}. \quad (8)$$

Based on stability theory, the traffic jam will never happen in this system if $p(s)$ is stable and $\|G(s)\|_\infty \leq 1$.

As $d^2 x_n(t)/dt^2$ is in positive correlation with $[V_{new}^{op}(\Delta x_n(t), v_n(t)) - v_n(t)]$, we have $\alpha > 0$. Because of the OV function characteristic, we have $\Lambda_1 > 0$, so the condition for $p(s)$ to be stable is $\lambda > \alpha\Lambda_2$.

Then we consider $\|G(s)\|_\infty \leq 1$ which can be expressed as

$$|G(j\omega)|^2 = |G(-j\omega)G(j\omega)|$$
$$= \frac{\alpha^2\Lambda_1^2 + \omega^2\lambda^2}{(\alpha\Lambda_1 - \omega^2)^2 + (\alpha+\lambda-\alpha\Lambda_2)^2\omega^2} \leq 1.$$

The sufficient condition can be obtained as

$$\alpha^2 + \alpha^2\Lambda_2^2 - 2\alpha^2\Lambda_2 + 2\alpha\lambda - 2\alpha\lambda\Lambda_2$$
$$-2\alpha\Lambda_1 + \omega^2 \geq 0, \quad (9)$$
$$\omega \in [0, \infty)$$

which can be rewritten as

$$\alpha + \alpha\Lambda_2^2 - 2\alpha\Lambda_2 + 2\lambda - 2\lambda\Lambda_2 - 2\Lambda_1 \geq 0. \quad (10)$$

The neutral stability surface in the space of $(\Lambda_1, \Lambda_2, \alpha)$ with different $b$ are shown in Fig.1. When $b=0$, the model equals the full velocity difference model (FVDM) [14], and the figure shows a neutral stability line in Fig.1(a). According to Fig.1(b) - (d), with higher $b$, the peak of stability surface moves to larger $\Delta x_n(t)$ when $v_n(t)$ increases. As a result, to avoid traffic jam at higher speed, driver needs to be more carefully to keep $\alpha$ above the stability surface.

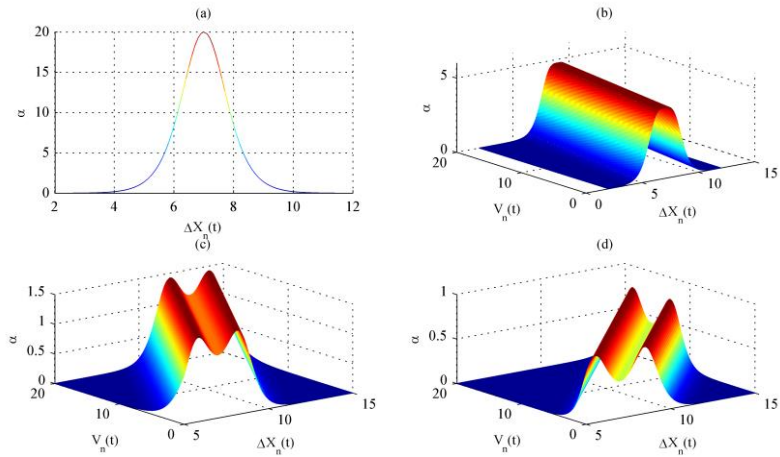

Fig.1 Neutral stability surface with different $b$: (a) $b=0$ (FVDM); (b) $b=0.1$; (c) $b=0.3$; (d) $b=0.5$.

## 4. Numerical simulations

Let us consider a 100-vehicle system running on an expressway without overtaking. In the simulations, the initial values are set as $v_n(0) = 15$m/s, $v_{max} = 20$m/s, $h_c = 7$m, $\alpha = \lambda = 0.5$.

Firstly, we assume in the simulation process the leading vehicle has a random fluctuation in acceleration or deceleration. Fig.2(a) shows the space-time plot when $b = 0$ (FVDM) and we can observe the oscillating headway distances in the lower vehicles. Fig.2(b) represents the new model with $b = 0.3$, indicating traffic jam never happens.

Secondly, we consider the situation with disturbance where the leading vehicle stops suddenly from t=110s to t=115s. The space-time plot and velocity behavior of three vehicles are described in Fig.3 - 5. We can see heavy traffic jams when $b = 0$, and as $b$ increases the traffic jams are suppressed.

Simulations illustrate that with a larger headway distance but higher speed, vehicles in new model are more sensitive to the variation of leading vehicle's acceleration or deceleration to keep the traffic system in stable situation.

## 5. Summary

In this paper, we proposed a new car-following model considering the variable safety headway distance based on the typical OVM. The stability condition of new model is analyzed by applying the control theory. Finally, numerical simulations are given, and the results are consistent with the theoretical study. In conclusion, the new model suppresses the traffic jam in typical OVM with high speed and larger headway distance.

## Acknowledgement

The research in this paper was conducted as part of the project " The occurrence and evolution of traffic gridlock in mega-city under storm rain conditions ", which is funded by Beijing Natural Science Foundation (Project No. 9132010).

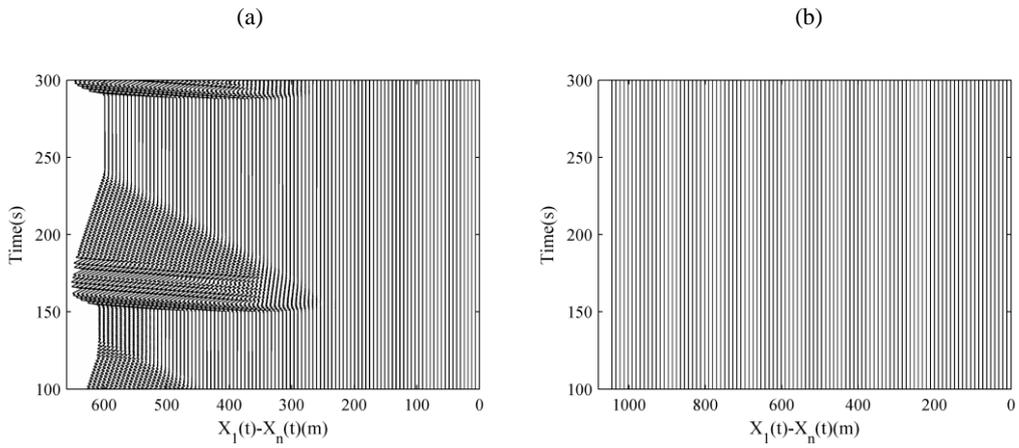

Fig.2 Space-time plot of: (a) Typical model (*b*=0); (b) New model (*b*=0.3).

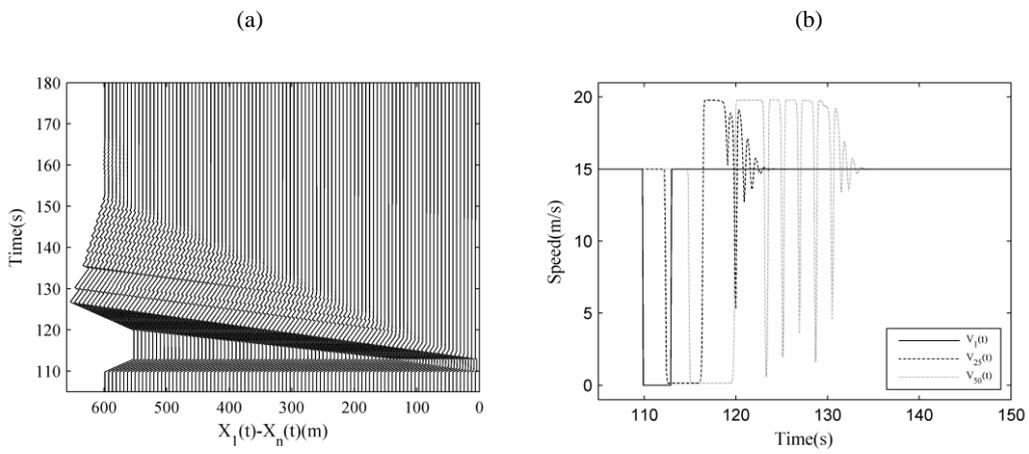

Fig.3 (a) Space-time plot; (b) Velocity of the first, 25$^{th}$ and 50$^{th}$ vehicles (*b*=0).

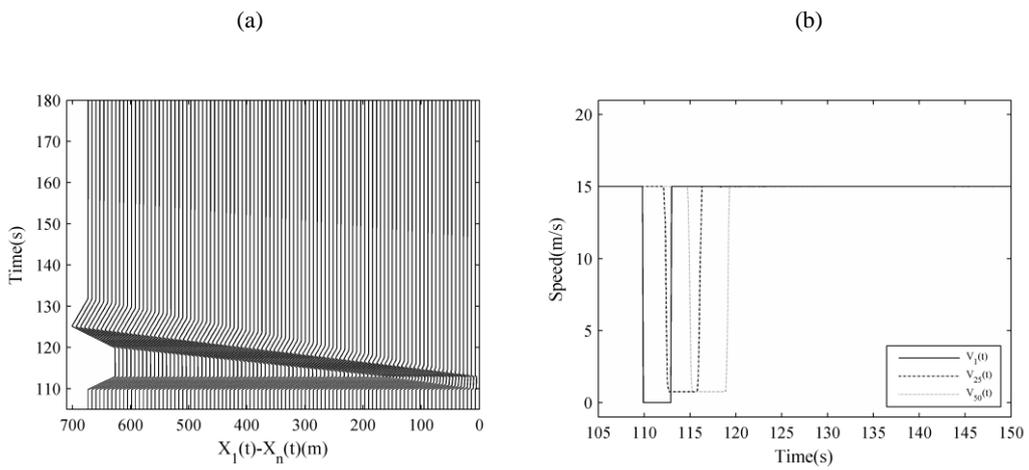

Fig.4 (a) Space-time plot; (b) Velocity of the first, 25$^{th}$ and 50$^{th}$ vehicles (*b*=0.05).

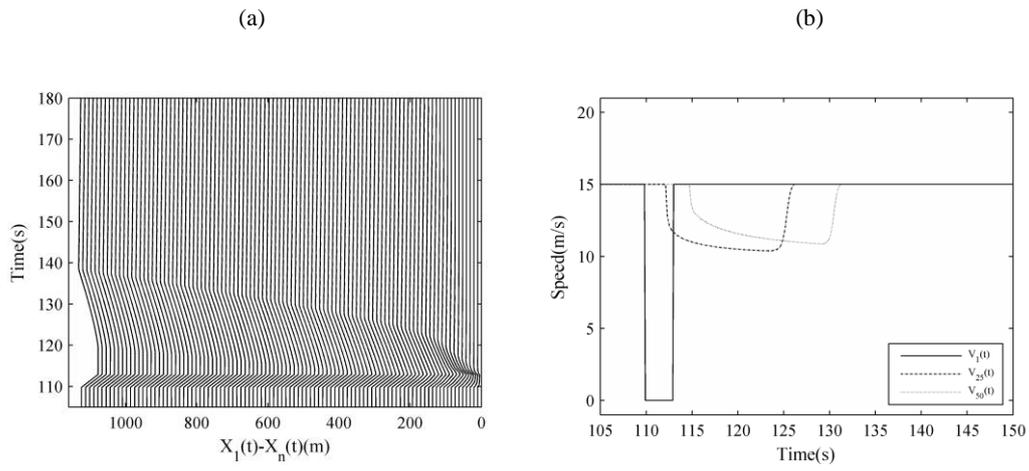

Fig.5 (a) Space-time plot; (b) Velocity of the first, 25$^{th}$ and 50$^{th}$ vehicles ($b$=0.3).

List of Figures Captions

Fig.1 Neutral stability surface with different *b*: (a) b=0 (FVDM); (b) *b*=0.1; (c) *b*=0.3; (d) *b*=0.5.

Fig.2 Space-time plot of (a) Typical model (*b*=0); (b) New model (*b*=0.3)

Fig.3 (a) Space-time plot (*b*=0); (b) Velocity of the first, $25^{th}$ and $50^{th}$ vehicles (*b*=0)

Fig.4 (a) Space-time plot (*b*=0.05); (b) Velocity of the first, $25^{th}$ and $50^{th}$ vehicles (*b*=0.05)

Fig.5 (a) Space-time plot (*b*=0.3); (b) Velocity of the first, $25^{th}$ and $50^{th}$ vehicles (*b*=0.3)